\documentclass[]{amsart}

\usepackage{amsfonts}
\usepackage{amscd}
\usepackage{epsfig}
\usepackage{epic,eepic}
\usepackage{graphicx}
\usepackage{psfrag}
\usepackage{amsmath,mathrsfs,amssymb}
\usepackage{amsthm}
\usepackage{setspace}
\usepackage{url}
\usepackage{here}
\usepackage{diagrams}

\newtheorem{theorem}{Theorem}[section]

\newtheorem{proposition}[theorem]{Proposition}
\newtheorem{corollary}[theorem]{Corollary}
\newtheorem{question}[theorem]{Question}

\theoremstyle{definition}
\newtheorem{remark}[theorem]{Remark}
\newtheorem{definition}[theorem]{Definition}

\newcommand{\PP}{\mathbb P^{1}}
\newcommand{\LL}{\mathscr L}
\newcommand{\OO}{\mathscr O}
\newcommand{\ZZ}{\mathbb Z}
\newcommand{\WW}{W^{r}_{d}(C)}

\newcommand{\AFG}{\mathcal A_{F}(G)}
\newcommand{\BWW}{\overline{W}^{r}_{2n}(C)} 
\newcommand{\BAA}{\overline{\mathcal A}_{F}(G)}
\newcommand{\MWW}{\mathcal W^{r}_{2n}}
\newcommand{\MWWO}{\mathcal W^{1}_{2n}}
\newcommand{\WWO}{W^{1}_{2n}(C)}
\newcommand{\WWN}{W^{r}_{2n}(C)}

\title{linear series on ribbons}

\author{Dawei Chen}

\address{University of Illinois at
  Chicago, Mathematics Department, Chicago, IL 60607}

\email{dwchen@math.uic.edu}


\begin{document}
\bibliographystyle{plain}

\begin{abstract}
A ribbon is a double structure on $\PP$. The geometry of a ribbon is closely related to that of a smooth curve. In this note we consider
linear series on ribbons. Our main result is an explicit determinantal description for the locus $W^{r}_{2n}$ of degree $2n$ line bundles with at least $(r+1)$-dimensional sections on a ribbon. We also discuss some results of Clifford and Brill-Noether type. 
\end{abstract}

\maketitle

\tableofcontents

\section{Introduction}
In this section, we recall some basic theory of ribbons. Many results and notations below come from Bayer-Eisenbud \cite{BE}. 

We work over an algebraically closed field $k$ of characteristic 0. 
A ribbon on $\PP$ is a scheme $C$ equipped with an isomorphism $\PP \rightarrow C_{red}, $ such that the ideal sheaf 
$\LL$ of $\PP$ in $C$ satisfies $$\LL^{2} = 0.$$ 
Because of this condition, $\LL$ can be regarded as a line bundle on $\PP$. It is called the \emph{conormal bundle}
of $\PP$ in $C$. There is a short exact sequence called the \emph{conormal sequence}: 
$$ 0\rightarrow \LL \rightarrow \OO_{C} \rightarrow \OO_{\PP} \rightarrow 0. $$ 

Define the arithmetic genus $g$ of a ribbon $C$ as: 
$$g = 1 - \chi(\OO_{C}). $$
From the conormal sequence, we see that $C$ has genus $g$ if and only if the conormal bundle $\LL$ on $\PP$ is isomorphic to $\OO_{\PP}(-g-1). $

There is another short exact sequence called the \emph{restricted cotangent sequence}: 
$$ 0 \rightarrow \LL \rightarrow \Omega_{C}|_{\PP} \rightarrow \Omega_{\PP} \rightarrow 0. $$ 
This restricted cotangent sequence defines the \emph{extension class} of $C$: 
$$e_{c} \in\  \mbox{Ext}^{1}_{\PP}(\Omega_{\PP}, \LL). $$ 

We will say that two ribbons are \emph{isomorphic} over $\PP$ if there is an isomorphism between them that extends the identity 
on $\PP$. A ribbon $C$ is \emph{split} if the inclusion $\PP \hookrightarrow C$ admits a section. Such a section is a scheme-theoretically degree two map from 
$C$ to $\PP$. We also call $C$ \emph{hyperelliptic} if it is split. 

\begin{theorem}
\label{main}
Given any line bundle $\LL$ on $\PP$ and any class $e \in Ext^{1}_{\PP}(\Omega_{\PP}, \LL)$, there is a unique ribbon 
$C$ on $\PP$ with $e_{c} = e$. If there is another class $e' \in Ext^{1}_{\PP}(\Omega_{\PP}, \LL)$ corresponding to a
ribbon $C'$, then $C\cong C'$ if and only if $e = ae'$ for some $a\in k^{*}$. A hyperelliptic ribbon corresponds to the split extension. 
The set of nonhyperelliptic ribbons of genus $g$, up to isomorphism over $\PP$, is the set 
$$\mathbb P^{g-3} = \mathbb P(H^{0}(\OO_{\PP}(g-3))). $$
\end{theorem}

\begin{proof}
The above results are essentially from \cite[Thm. 1.2, 2.1]{BE}. 
\end{proof}

There is an explicit way to write down the structure of a ribbon by a gluing method, cf. \cite[Sec. 3]{BE}. 
Define two open sets 
$$u_{1} = \mbox{Spec}\  k[s], \ u_{2} = \mbox{Spec}\  k[t]$$
that cover $\PP$ via the identification $s^{-1} = t$ on $u_{1}\cap u_{2}.$ 

If $C$ is a genus $g$ ribbon on $\PP$, we can write 
$$U_{1}: = C|_{u_{1}} \cong \mbox{Spec}\  k[s, \epsilon]/\epsilon^{2}, $$
$$U_{2}: = C|_{u_{2}} \cong \mbox{Spec}\  k[t, \eta]/\eta^{2}. $$
$C$ may be specified by a gluing isomorphism between $U_{1}$ and $U_{2}$ over $u_{1} \cap u_{2}.$
The ideal sheaf $\LL\cong \OO_{\PP}(-g-1)$ of $\PP$ in $C$ is generated by $\epsilon$ on $u_{1}$ and by $\eta$ on $u_{2}$. 
So we can further write 
$$ \epsilon = t^{-g-1} \eta, $$
$$ s^{-1} = t + F(t)\eta$$ 
on $u_{1}\cap u_{2}$, with $F(t) \in k[t, t^{-1}].$
Conversely, any such gluing data can determine a ribbon of genus $g$ on $\PP$. 

If we change the coordinates with 
$$s' = s + p(s)\epsilon, \ t = t'+q(t)\eta $$
on $U_{1}$ and $U_{2}$ with polynomials $p(s), q(t)$,  
then we get 
\begin{eqnarray*}
s'^{-1} & =  & s^{-1} - s^{-2}p(s)\epsilon   \\
            & = & t' + \big( F(t) - t^{-g+1}p(t^{-1}-F(t)\eta)+q(t) \big)\eta  \\
            & = & t + F(t)\eta - (t+F(t)\eta)^{2} p(t^{-1}) t^{-g-1}\eta \\
            & = & t + F(t)\eta - t^{1-g}p(t^{-1}) \eta \\ 
            & = & t' + \big(F(t) + q(t) - t^{1-g}p(t^{-1})\big)\eta \\
            & = & t' + \big(F(t') + q(t') - t'^{1-g}p(t'^{-1})\big)\eta
\end{eqnarray*}
The fact that $t\eta = t'\eta$ is used in the last step. 
If we multiply $s$ or $t$ by a scalar, $F$ will also be multiplied by the same scalar. 
Therefore, $F$ can be determined as an element of the projective space of lines in the quotient 
$$k[t, t^{-1}] / (k[t] + t^{-g+1}k[t^{-1}]). $$
From now on, we shall write $F$ as 
\begin{eqnarray}
\label{F}
F & = & \sum_{i=1}^{g-2}F_{i} t^{-i}. 
\end{eqnarray}
$F = 0$ corresponds to a hyperelliptic ribbon. This explicit expression recovers the fact in Theorem~\ref{main} that 
nonhyperelliptic ribbons of genus $g$, up to isomorphism over $\PP$, is parameterized by the set 
$$\mathbb P^{g-3} = \mathbb P(H^{0}(\OO_{\PP}(g-3))). $$ 

Now we consider line bundles. Let $L$ be a line bundle on a ribbon $C$. The degree of $L$ 
is defined as 
$$\mbox{deg}\  L: = \chi(L) - \chi (\OO_{C}). $$ 

\begin{proposition}
\label{line bundle}
 If $L|_{\PP} \cong \OO_{\PP}(n)$, then deg $L$ = $2n$. 
The Picard group of $C$ 
$$ \mbox{Pic}\  C = H^{1}(\OO_{\PP}(-g-1)) \times \ZZ \cong k^{g} \times \ZZ, $$
where the projection to $\ZZ$ is given by the degree of the restriction $L|_{\PP}. $
\end{proposition}

\begin{proof}
See \cite[Prop. 4.1, 4.2]{BE}. 
\end{proof}

Bayer and Eisenbud remarked that one must switch to torsion-free sheaves in order to obtain the analogue of line bundles of odd degree. For simplicity, 
here we only consider line bundles of even degree $2n$. 

A line bundle $L$ can also be constructed by gluing. Suppose $L|_{\PP} \cong \OO_{\PP}(n)$. 
Using the above notation, we have 
$$ L|_{U_{1}} = k[s,\epsilon]e_{1}, $$
$$ L|_{U_{2}} = k[t,\eta]e_{2}, $$
and 
$$ e_{1} = (t+F\eta)^{n}(1+G\eta)e_{2}$$ 
on $U_{1}\cap U_{2}$ for some $G \in k[t , t^{-1}]. $
Conversely, any such $G$ can be used to construct a line bundle. 

If we change the coordinates by
$$ e_{1} = (1+m(s)\epsilon)^{-1} e'_{1}, \ e_{2} = (1+n(t)\eta) e'_{2} $$ 
on $U_{1}$ and $U_{2}$ with polynomials $m(s), n(t)$, then we get 
$$e'_{1} = (t+F\eta)^{n} \big(1+(G-m(t^{-1})t^{-g-1} + n(t))\eta\big) e'_{2}. $$
In order to classify $L$, it suffices to specify $G$ as an element of 
$$ k[t, t^{-1}] / (k[t] + t^{-g-1}k[t^{-1}]) = H^{1}(\OO_{\PP}(-g-1)). $$
This also recovers the fact in Proposition~\ref{line bundle} that $H^{1}(\OO_{\PP}(-g-1))$ parameterizes line bundles of fixed degree. 
We will also write $G$ as 
\begin{eqnarray}
\label{G}
G & = & \sum_{j = 1}^{g} G_{j}t^{-j}. 
\end{eqnarray}
The line bundle corresponding to $G=0$ is isomorphic to the pullback of $\OO_{\PP}(n)$ from $C_{red}\cong \PP$ to $C$. 

Let $L$ be a line bundle on $C$ of degree $2n$ given by the above gluing data. 
We would like to find out the space of global sections of $L$. 
Suppose $p = p(t)$ is a polynomial of degree $\leq n$. $pe_{2}|_{\PP}$ determines
an element 
$$ \sigma \in H^{0}(L|_{\PP}) = H^{0}(\OO_{\PP}(n)). $$
Define 
\begin{eqnarray}
\label{delta}
\delta_{L}(p) & = & -(p'F + pG) \in k[t,t^{-1}]/(k[t] + t^{n-g-1}k[t^{-1}]), 
\end{eqnarray}
where $p' = \frac{\partial p}{\partial t}. $

\begin{theorem}
\label{section}
The space of sections of $L$ restricted to $U_{2} = Spec\  k[t,\eta]$ can be identified as the direct sum
of the space of elements $q(t)\eta$ and the space of expressions $p(t) + p_{1}(t)\eta$, 
where $q$ is a polynomial of degree $\leq n-g-1$, $p$ is a polynomial of degree $\leq n$ 
satisfying $\delta_{L}(p) = 0$ in (\ref{delta}) and $p_{1}\in k[t]$ is the polynomial part of $p'F + pG$, 
i.e. $p_{1}(t) \equiv p'(t) F(t) + p(t) G(t)$ mod $t^{-1}k[t^{-1}].$ 
\end{theorem}

\begin{proof}
This is exactly \cite[Thm. 4.3]{BE}.
\end{proof}

At a first glance, the above way to identify $H^{0}(L)$ seems quite messy. Nevertheless, a further observation will imply an important 
conclusion immediately. 

\begin{corollary} 
\label{2g}
Let $L$ be a line bundle of degree $2n$ on a ribbon $C$. If $n \geq g$, then $h^{0}(L) = 1- g + 2n$. 
\end{corollary}

\begin{proof} Let a section of $L$ restricted to $U_{2}$ correspond to the data $(q(t)\eta, p(t) + p_{1}(t))$ in the above setup. 
When $n\geq g$, we have $k[t,t^{-1}]/(k[t] + t^{n-g-1}k[t^{-1}]) \equiv 0$. Then by its definition, 
$\delta_{L}$ always takes value on 0. So $\delta_{L}(p) = 0$ does not impose any condition on $p(t)$. 
The only constraint of $p(t)$ and $q(t)$ is the upper bound of their degree. $q(t)$ has degree $\leq n-g -1$ and 
$p(t)$ has degree $\leq n$. In total, they yield a $(1-g+2n)$-dimensional space for sections of $L$. 
\end{proof}

\begin{remark}
Notice that if $n\geq g$, then the degree $d$ of $L$ satisfies $d=2n>2g-2$. In case $C$ is a smooth curve of genus $g$, 
$h^{0}(L) = 1-g+d$ holds for any line bundle $L$ on $C$ with degree $d>2g-2$. So the above corollary can be viewed as a similarity 
between ribbons and smooth curves. 
\end{remark}

\section{The locus $W^{r}_{2n}$}

For smooth curves, the theory of special linear series can be best characterized by the Brill-Noether theory. In order to avoid missing any 
work, we simply refer readers to \cite[Chap. V]{ACGH} for bibliographical notes on this topic. Let $C$ be a curve of genus $g$. We introduce the variety 
$\WW$ as 
$$\WW = \{L\in \mbox{Pic}^{d}(C)\ |\  h^{0}(C,L) \geq r+1   \}. $$
We also define the Brill-Noether number $\rho$ as
$$\rho = g - (r+1)(g-d+r).$$
The basic results of the Brill-Noether theory can be summarized as follows. 

\begin{theorem}
\label{BN}
 Let $C$ be a smooth curve of genus $g$. Let $d, r$ be integers such that $d\geq 1, r\geq 0$.

(Existence). If
$\rho \geq 0$, $\WW$ is non-empty. Furthermore, every component of $\WW$ has dimension at least equal to $\rho$ provided $r \geq d-g$. 

(Connectedness). Assume that $\rho \geq 1$. Then $\WW$ is connected. 

(Dimension). Assume that $C$ is a general curve. If $\rho < 0$, then $\WW$ is empty. 
If $\rho \geq 0$, then $\WW$ is reduced and of dimension $\rho$. 
\end{theorem}

We would like to investigate linear series for ribbons. The importance of such a study has three folds. In the first place, $\WW$ essentially carries
a determinantal structure for a smooth curve $C$. In case $C$ is a ribbon, the determinantal characterization can even be made explicit. Secondly, 
the Brill-Noether theory for a special member in a family of curves usually reveals information for a general one. Ribbons do arise as the degeneration of smooth curves, cf. \cite{F}. Finally, Lazarsfeld \cite{L} proved that a general curve contained in certain K3 surface satisfies the above dimension theorem. 
Correspondingly, ribbons lie on the so-called K3 carpet, i.e. double structure on a rational normal scroll, which has the same numerical invariants as a 
smooth K3 surface, cf. \cite[Sec. 8]{BE}. Hence, it would be interesting to figure out some results of Brill-Noether type for ribbons. 

Let $C$ be a ribbon determined by $[F_{1}, \ldots, F_{g-2}]$ the coefficients of $F$ in (\ref{F}) up to scalar. Let $L$ be a line bundle of degree $2n$ on $C$ 
determined by $(G_{1}, \ldots, G_{g})$ the coefficients of $G$ in (\ref{G}). If $n \geq g$, there is no special linear system because of Corollary~\ref{2g}. 
Actually we only need to consider $2n \leq g - 1$, since the Riemann-Roch formula also holds for ribbons, cf. \cite[Sec. 5]{BE}. From now on, assume that $2n\leq g-1$. 
Define a $(g-n)\times (n+1)$ matrix $\AFG$ with entries $s_{ij} = G_{i+j-1} + (j-1)F_{i+j-2}$, namely, 
\[ \AFG = \left( \begin{array}{cccc}
G_{1} & G_{2} + F_{1} & \cdots & G_{n+1}+nF_{n} \\
G_{2} & G_{3} + F_{2} & \cdots & G_{n+2} + nF_{n+1} \\
\vdots & \vdots              & \ddots &  \vdots        \\ 
G_{g-n} & G_{g-n+1} + F_{g-n} & \cdots & G_{g} 
  \end{array} \right)  \] 
Now we can state our main result. 

\begin{theorem} 
\label{det}
In the above setting, the locus $\WWN$ is isomorphic to the following affine algebraic set: 
$$ \WWN = \{ (G_{1}, \ldots, G_{g}) \in \mathbb A^{g}\  |\  \mbox{rank}\  \AFG  \leq n - r \}. $$ 
\end{theorem}

\begin{proof}
By Theorem~\ref{section}, the space $H^{0}(L)$ of global sections of $L$ can be identified as the direct sum of two spaces: 
$$ \langle q(t)\eta \rangle \oplus \langle p(t) + p_{1}(t)\eta \rangle. $$
Since $q(t)$ is a polynomial of degree $\leq n - g - 1$ and $n \leq g-2$, the first summand is the null space. For the second, $p(t)$ is a polynomial of degree 
$\leq n$ satisfying $\delta_{L}(p) = 0$ as in (\ref{delta}). Then $p_{1}(t)$ is determined by $p(t)$ the polynomial part of $p'F + pG$. 
Let $p(t) = \sum_{i = 0}^{n} a_{i} t^{i}.$ The condition $\delta_{L}(p) = 0$ means 
$$p'F+pG \in k[t] + t^{n-g-2}k[t^{-1}],$$ which is equivalent to the following:
$$\AFG \cdotp \vec{a} = 0, $$
where $\vec{a}$ is the vector $(a_{0}, \ldots, a_{n})^{t}$ determined by the coefficients 
of $p(t)$. Note that $g-n \geq n+1$ by the assumption on $n$. Hence, $\WWN$ can be identified as the desired determinantal locus. 
\end{proof}
 
The following Clifford theorem for ribbons is a direct consequence of Theorem~\ref{det}. 
 
\begin{theorem}
Let $C$ be a ribbon and $L$ be a line bundle of degree $2n$ on $C$, $1\leq n \leq g-2$. Then $h^{0}(C, L) \leq n+1$. The equality holds 
if and only if $C$ is a hyperelliptic ribbon and $L$ is the pullback of $\OO_{\PP}(n)$ from $C_{red}\cong \PP$ to $C$. 
\end{theorem}

\begin{proof}
By the above determinantal description for $\WWN$, we know that $h^{0}(C, L) \leq n+1$. If the equality holds, then $r = n$. We have 
$G_{i} = 0$ and $F_{j} = 0$ for all $i, j$. Thus $C$ is hyperelliptic and $L$ is isomorphic to the pullback of $\OO_{\PP}(n)$ from $C_{red}\cong \PP$. 
\end{proof}

When $C$ is a hyperelliptic ribbon, i.e. $F_{i} = 0$ for all $i$, $\AFG$ has entries $s_{ij} = G_{i+j-1}$: 
\[ \left( \begin{array}{cccc}
G_{1} & G_{2} & \cdots & G_{n+1} \\
G_{2} & G_{3} & \cdots & G_{n+2} \\
\vdots &\vdots & \ddots &  \vdots        \\ 
G_{g-n} & G_{g-n+1} & \cdots & G_{g} 
  \end{array} \right) \]
  
Such a matrix is called the Catalecticant matrix. 
We cite the result \cite[Prop. 4.3]{E} as follows. 
\begin{proposition}
\label{cat} 
The space of rank $\leq m$ Catalecticant matrices is isomorphic to a cone over $S_{m}$, where $S_{m}$ is the union of $m$-secant $(m-1)$-planes 
to a rational normal curve of degree $g-1$. 
\end{proposition}

This exactly describes $\WWN$ for a hyperelliptic ribbon. 

\begin{theorem}
\label{hyperelliptic}
Let $C$ be a hyperelliptic ribbon. Then $\WWN$ is isomorphic to a cone over $S_{n-r}$ for $r < n$. In particular, $\WWN$ has dimension equal to $2n-2r$. 
\end{theorem}

\begin{proof}
$\WWN$ can be identified as the space of rank $\leq n-r$ Catalecticant matrices, which is isomorphic to a cone over $S_{n-r}$ by Proposition~\ref{cat}. 
$S_{n-r}$ has dimension $2n-2r-1$, so $\WWN$ has dimension $2n-2r$. 
\end{proof}

We have seen that for a nonhyperelliptic ribbon, its structure can be determined by the data $[F_{1}, \ldots, F_{g-2}]$ in (\ref{F}) as a point of $\mathbb P^{g-3}$. 
Note that the expected dimension of $\WWN$ would still be $g - (r+1)(g-2n+r)$, which equals the Brill-Noether number $\rho$. We would like to study the actual dimension 
of $\WWN$. Firstly, let us focus on a natural 
compactification of $\WWN$ as follows. 

Define another $(g-n)\times (n+1)$ matrix $\BAA$ with entries $s_{ij} = G_{i+j-1} + (j-1)F_{i+j-2}G_{0}$:
\[ \BAA = \left( \begin{array}{cccc}
G_{1} & G_{2} + F_{1}G_{0} & \cdots & G_{n+1}+nF_{n}G_{0} \\
G_{2} & G_{3} + F_{2}G_{0} & \cdots & G_{n+2} + nF_{n+1}G_{0} \\
\vdots & \vdots              & \ddots &  \vdots        \\ 
G_{g-n} & G_{g-n+1} + F_{g-n}G_{0} & \cdots & G_{g} 
  \end{array} \right)  \] 
 
Let $$\BWW =  \{ [G_{0}, G_{1}, \ldots, G_{g}] \in \mathbb P^{g} \ | \ \mbox{rank}\  \BAA  \leq n - r \}. $$  
There is an inclusion $\WWN \subset \BWW$ given by 
$$ (G_{1}, \ldots, G_{g}) \rightarrow [1, G_{1}, \ldots, G_{g}] . $$
The complement of $\WWN$ in $\BWW$ is just the hyperplane section $\{G_{0} = 0\} \cap \BWW$. 

We also need to introduce generic determinantal varieties. Let $M$ be the space of $(g-n)\times (n+1)$ matrices. Denote by $M_{l}$ the locus of rank $\leq l$ matrices. $M_{l}$ is called the $l$-generic determinantal variety, $l \leq n+1$. Denote by $\mathbf M$ and $\mathbf M_{l}$ the projectivization of 
$M$ and $M_{l}$ respectively. 

\begin{proposition}
$\mathbf M_{l}$ is an irreducible subvariety of codimension $(g-n-l)(n+1-l)$ in $\mathbf M$. 
\end{proposition}

One can refer to \cite[Chap. II]{ACGH} for a general discussion on determinantal varieties. Our next result is about the global geometry of 
$\BWW$. 

\begin{theorem} 
\label{bww}
Let $C$ be a ribbon of genus $g$. For $r < n$, $\BWW$ is always non-empty and has dimension equal to $2n-2r-1$ or $2n-2r$. Each irreducible component of $\BWW$ has dimension at least equal to $\rho$ provided $\rho \geq 0$. 
Furthermore,  $\BWW$ is connected provided $\rho > 0$. 
\end{theorem}

\begin{proof} $\BWW$ is the intersection of $\mathbf M_{n-r}$ and a $g$-dimensional linear subspace of $\mathbf M$ determined by $s_{ij} = G_{i+j-1} + (j-1)F_{i+j-2}G_{0}$. Therefore, each irreducible component of $\BWW$ has dimension $\geq$ dim $\mathbf M_{n-r} + g - $ dim $\mathbf M = \rho $. 

When $G_{0} = 0$, the matrix $\BAA$ reduces to a Catalecticant matrix with entries $s_{ij} = G_{i+j-1}$. The space of rank $\leq n-r$ Catalecticant matrices
has dimension $2n-2r$ by Proposition~\ref{cat}. It implies that the hyperplane section $\{G_{0} = 0\} \cap \BWW$ has dimension $2n-2r-1$. If the top dimensional 
component of $\BWW$ is contained in $\{G_{0} = 0\}$, then $\BWW$ has dimension $2n-2r-1$. Otherwise it has dimension $2n-2r$. 

$\BWW$ can also be regarded as the intersection of the $(n-r)$-generic determinantal variety $\mathbf M_{n-r}$ and a $g$-dimensional linear subspace 
of $\mathbf M$ defined by $s_{ij} =  G_{i+j-1} + (j-1)F_{i+j-2}G_{0}$ for a fixed lifting $(F_{1}, \ldots, F_{g-2})$. If $\rho > 0$, the sum of the dimensions of these two spaces is greater than the dimension of $\mathbf M$. The connectedness of their intersection $\BWW$ follows as a consequence of \cite[Ex. 3.3.7]{L}. 
\end{proof}

\begin{corollary}
\label{at least}
Assume that $\rho \geq 0$. If $\WWN$ exists for a ribbon $C$, then $\WWN$ has dimension at least equal to $\rho$. 
\end{corollary}

\begin{proof}
$\WWN$ is the complement of the hyperplane section $\{G_{0} = 0\} \cap \BWW$ in $\BWW$. By Theorem~\ref{bww} we know that each component of 
$\BWW$ has dimension $\geq \rho$. So does $\WWN$ assuming it is non-empty.
\end{proof}

We can also let $(F_{1}, \ldots, F_{g-2})$ vary as a point of $\mathbb A^{g-2}$ and define the global Brill-Noether locus $\MWW$ 
as follows: 
\[ \MWW = \{(G_{1}, \ldots, G_{g}; F_{1}, \ldots, F_{g-2}) \in \mathbb A^{g} \times \mathbb A^{g-2} \ |\ \mbox{rank}\ \AFG \leq n-r \}. \]
$\MWW$ is the intersection of the $(n-r)$-generic determinantal variety $M_{n-r}$ and a $(2g-2)$-dimensional linear subspace $S$ of $M$, where $S$ is determined by 
relations $2s_{ij} = s_{i-1\ j+1} + s_{i+1\ j-1}.$ Note that the expected dimension of $\MWW$ would be $2g-2-(g-2n+r)(r+1) = g-2 + \rho$, which implies the following
conclusion right away.

\begin{corollary}
\label{fiber}
If $\MWW$ has dimension equal to $g-2+\rho$, then for $(F_{1}, \ldots, F_{g-2})$ corresponding to a general ribbon $C$, $\WWN$ has dimension at most equal to $\rho$. 
\end{corollary}

In order to calculate the actual dimension of $\MWW$, we introduce the concept of $l$-generic spaces developed by Eisenbud \cite[Prop.-Def. 1.1]{E}. 

\begin{definition}
Let $L$ be a linear subspace of the space $M$ of $(g-n)\times (n+1)$ matrices. $L$ can be regarded as an associated $(g-n)\times (n+1)$ matrix of linear forms. We say that $L$ is $m$-\emph{generic} for some integer $1\leq m \leq n+1$ if after arbitrary invertible row and column operations, any $m$ of the linear forms $L_{ij}$ in $L$ are linear independent. 
\end{definition}

We also say that $L$ meets $M_{l}$ \emph{properly} if their intersection has codimension equal to $(g-n-l)(n+1-l)$ in $L$. 

\begin{theorem}
\label{proper}
Let $L\subset M$ be a $m$-generic space, then $L$ meets $M_{n+1-m}$ properly. 
\end{theorem}

\begin{proof}
This is part of \cite[Thm. 2.1]{E}. 
\end{proof}

Note that the space of Catalecticant matrices is $1$-generic. One can also prove the 2-genericity for the space $S$ of matrices of type $\AFG$. 

\begin{proposition}
Consider $G_{1}, \ldots, G_{g}; F_{1}, \ldots, F_{g-2}$ as independent linear forms. The (2g-2)-dimensional vector space $S$ of all matrices determined by $\AFG$ is 2-generic.
\end{proposition}

\begin{proof}
$\AFG$ is the matrix with entries $s_{ij} = G_{i+j-1} + (j-1)F_{i+j-2}$.  
Suppose there were two invertible matrices $A = (a_{ij})$ and $B = (b_{ij})$ corresponding to invertible row and column operations such that  
two entries of the new matrix $A\cdot\AFG\cdot B = (s'_{ij})$ became equal to each other. We can always assume that these two entries are $s'_{11}$ and $s'_{22}$. The case that they are in the same row or column would be even easier. Then the condition $s'_{11} = s'_{22}$ is equivalent to
$$\sum_{i+j = k+1}(a_{1i}b_{j1}-a_{2i}b_{j2})\big(G_{k}+(j-1)F_{k-1}\big) = 0 $$
for any $k$. Namely, 
$$\sum_{i+j=k+1}(a_{1i}b_{j1}-a_{2i}b_{j2}) = 0 \ \mbox{and}\ \sum_{i+j=k+1}(a_{1i}b_{j1}-a_{2i}b_{j2})(j-1) = 0$$
since $G_{i}$ and $F_{j}$ can vary independently. 

Define four polynomials as follows: 
$$A_{k}(x) = \sum_{i}a_{ki}x^{i} \ \mbox{and}\ B_{k}(x) = \sum_{j}b_{kj}x^{j} \ \mbox{for}\ k=1,2. $$
We can deduce from the above two equalities that 
$$ A_{1}(x)B_{1}(x) = A_{2}(x)B_{2}(x) \ \mbox{and}\ A_{1}(x)B'_{1}(x) = A_{2}(x)B'_{2}(x). $$
Since the matrices $A$ and $B$ are invertible, by these two relations we can get 
$$ B_{1}(x)B'_{2}(x) = B_{2}(x)B'_{1}(x), $$ 
which would imply that $(B_{1}(x)/B_{2}(x))' = 0$. Then $B_{1}(x)/B_{2}(x)$ would be a constant, which contradicts to the assumption that 
the matrix $B$ is invertible. 
\end{proof}

\begin{corollary}
\label{at most}
For $r = 1$, $\MWWO$ has dimension $4n-4$ and $\WWO$ has dimension at most equal to $\rho = 4n - g - 2$ for a general ribbon $C$, provided $\rho \geq 0$. 
\end{corollary}

\begin{proof}
Since the space of matrices $\BAA$ is $2$-generic, it intersects $M_{n-1}$ properly. So the intersection $\MWWO$ has dimension equal to 
$g-2 + \rho = 4n-4$ by Theorem~\ref{proper}. Then by Corollary~\ref{fiber}, $\WWO$ has dimension at most equal to $\rho = 4n - 2- g$. 
\end{proof}

Based on Corollary~\ref{at least} and~\ref{at most}, we obtain the following conclusion. 

\begin{corollary}
For $r = 1$, if $\WWO$ exists for a general ribbon $C$, then $\WWO$ has dimension equal to $\rho = 4n-2-g$ provided $\rho \geq 0$. 
\end{corollary}

It would be interesting to pin down the following question in general.

\begin{question}
For $\rho \geq 0$, is the dimension of $\MWW$ equal to the expected dimension $g-2 + \rho$? For a general ribbon $C$, 
is the locus $\WWN$ non-empty and has dimension equal to $\rho$ provided $\rho \geq 0$? 
\end{question}

By the determinantal descriptions for $\MWW$ and $\WWN$, using Macaulay one can check that the above question
does have a positive answer when the genus of $C$ is small. 

{\bf Acknowledgements.} I am grateful to David Eisenbud for enlightening discussions on determinantal varieties and to Gregory Smith for his help on Macaulay. I also want to thank MSRI providing a wonderful environment for me to finish this paper.

\end{document}